\documentclass[pdflatex,sn-mathphys-num]{sn-jnl}

\usepackage{graphicx}%
\usepackage{multirow}%
\usepackage{amsmath,amssymb,amsfonts}%
\usepackage{amsthm}%
\newtheorem{theorem}{Theorem}[section] 
\newtheorem*{theorem*}{Theorem}
\usepackage{mathrsfs}%
\usepackage[title]{appendix}%
\usepackage{xcolor}%
\usepackage{textcomp}%
\usepackage{manyfoot}%
\usepackage{booktabs}%
\usepackage{algorithm}%
\usepackage{algorithmicx}%
\usepackage{algpseudocode}%
\usepackage{listings}%

\theoremstyle{thmstyleone}%
\newtheorem{proposition}[theorem]{Proposition}%

\theoremstyle{thmstyletwo}%

\theoremstyle{thmstylethree}%
\newtheorem{definition}{Definition}%

\theoremstyle{thmstylethree}%
\newtheorem{lemma}{Lemma}%

\theoremstyle{thmstylethree}%
\newtheorem{notation}{Notation}%

\begin{document}

\title[Article Title]{On graphical partitions with restricted parts}


\author{\fnm{Gilad} \sur{Levy}}


\abstract{An integer partition of $n$ is called graphical if its parts form a degree sequence of a simple graph. While unrestricted graphical partitions have been extensively studied, much less is known when the parts are restricted to a prescribed set. In this work, we investigate the probability that a uniformly random partition of an even integer $n$, subject to such restrictions, is graphical. We establish an upper bound on this probability expressed solely in terms of the Durfee square of the partition. Additionally, we evaluate the most-probable size of the Durfee square, depending solely on $n$ and the restrictions on the parts. Thus, we establish typical bounds on the probability that such random partitions are graphical. Our approach employs the Nash–Williams graphical condition, the saddle-point method and Edgeworth expansions.}

\keywords{graphical partitions, restricted parts, Durfee square, partitions.}



\maketitle

\section{Introduction}\label{sec1}
A graphical partition is an integer partition whose parts represent a degree sequence of a simple graph. This article studies graphical partitions with parts restricted to prescribed sets. We let $\mu(i)\in\mathbb{N}$ indicate the $i$-th smallest part a partition can have under the restrictions imposed on it, for all $i\in\mathbb{N}$ and some function $\mu:\mathbb{R}\to\mathbb{R}$. In this study, we only consider functions $\mu$ from the following set $M$.
\begin{definition}(set $M$)\label{def1}
     Let $M$ be the set of continuous, differentiable and strictly increasing functions $\mu:\mathbb{R}\to\mathbb{R}$ such that $\mu(0)=0$, $\mu(i)\in\mathbb{N}, \ \forall i\in\mathbb{N}$. 
\end{definition}
Representing discrete restrictions using continuous functions will allow us to apply analytic techniques. Additionally, For each $\mu \in M$, we define a corresponding set of partitions characterized by the restrictions imposed on their parts.
\begin{definition}(set $\mathcal{P}_\mu$)
    Let $\mu\in M$, then $\mathcal{P}_\mu$ is defined as the set of partitions with parts restricted to $\{\mu(i):i\in\mathbb{N}\}$.
\end{definition}

\newpage

Notice that $\mathcal{P}_\mu$ is a set of integer partitions whose $i$-th smallest possible part is $\mu(i)$, for all $i\in\mathbb{N}$, and for every $\mu\in M$. \\
\\
In our work, we find an upper bound for the probability that a random restricted partition is graphical. A key feature of our approach is that the resulting bound depends solely on the side length $D$ of the Durfee square of the partition and is invariant of the restrictions placed on the parts. This is formalized in the following theorem.
\begin{theorem*}[\ref{theo1}]\label{theo1}
    Let $\mu\in M$ and let $n\in\mathbb{N}$ be an even integer. Consider a partition of $n$ with Durfee square side length $D$, chosen uniformly at random from $\mathcal{P}_\mu$. Then, the probability that this partition is graphical is bounded above by
    \begin{equation}
        \frac{C}{\log^{3/2}D}
    \end{equation}
    for sufficiently large $D$, with $C=\sqrt{\frac{32}{\pi}}e^{-5/4}\approx 0.9144$.
\end{theorem*}
Theorem~\ref{theo2} establishes a trancscendental equation for the typical size of the Durfee square side length, depending directly on $n$ and the restrictions on the parts. \\
Using both theorems, we provide typical bounds on the probability that a random restricted partition will be graphical. We highlight the following special cases:
\begin{theorem*}[\ref{theo3}]\label{theo3}
    The probability that a random partition of $n$ is graphical, given that its parts are taken from a set whose elements grow asymptotically linearly, is bounded from above by
\begin{equation}
    \frac{C_\mathrm{lin}}{\log^{3/2}n}\Big(1+O\Big(\frac{1}{\log n}\Big)\Big)
\end{equation}
where $C_\mathrm{lin}=\sqrt{\frac{8}{\pi}}e^{-5/4}\approx0.4572$.
\end{theorem*}

\begin{theorem*}[\ref{theo4}]\label{theo4}
    The probability of a random partition of $n$ with parts restricted to perfect squares to be graphical is bounded from above by
    \begin{equation}
    \frac{C_\mathrm{sq}}{\log^{3/2}n}\Big(1+O\Big(\frac{\log\log n}{\log n}\Big)\Big)
\end{equation}
where $C_\mathrm{sq}=\sqrt{\frac{729}{2\pi}}e^{-5/4}\approx3.086$.
\end{theorem*}
\vspace{0.4cm}
The remainder of this paper is organized as follows. In Section 2, we provide the necessary preliminaries. Section 3 is dedicated to several probabilistic lemmas and estimates that underpin our main arguments. Finally, in Section 4, we provide the proof of our main result, Theorem~\ref{theo1}, and establish Theorem~\ref{theo2} to highlight important special cases in Theorems~\ref{theo3} and~\ref{theo4}

\newpage

\section{Preliminaries}
We recall the definition of successive ranks of an integer partition, as introduced by Atkin [1].
\begin{definition}[successive ranks]
    Let $k\in\mathbb{N}$ and let $\lambda$ be an integer partition. Denote by $X_k$ the number of parts of $\lambda$ that are at least $k$, and denote by $Y_k$ the $k$-th largest part in the partition. Then the $k$-th successive rank of $\lambda$ is defined by the difference $R_k=Y_k-X_k$.
\end{definition}
Furthermore, note the following definition of the Durfee square of an integer partition [2].
\begin{definition}[Durfee square]
    The Durfee square of an integer partition is the largest square that can fit in its Ferrers Diagram.
\end{definition}
In this study, we use the Nash-Williams condition for graphical sequences [3], which was later reformulated by Barnes and Savage [4], to connect between Number Theory and Graph Theory. 
\begin{theorem}[Nash-Williams]\label{Nash}
    An integer partition with Durfee square side length $D\in\mathbb{N}$ is graphical if and only if its successive ranks satisfy for all $1\le d\le D$,
    \begin{equation*}
        \sum_{k=1}^dR_k\le-d.
    \end{equation*}
\end{theorem}
Lastly, we recall the following theorem regarding Edgeworth expansions [5].
\begin{theorem}[Edgeworth expansion]\label{Edge}
    Let $d\in\mathbb{N}$, and let $A_k$ be independent random variables for all integers $1\le k\le d$, such that
\begin{equation*}
    \mathbb{E}[A_k]=0, \ \mathbb{E}[A_k^2]=\sigma_k^2, \ \mathbb{E}[A_k^3]<\infty.
\end{equation*}
Denote by $\kappa_{j,k}$ the $j$-th comulant of $A_k$, and let $\lambda_j=\sum_{k=1}^d\kappa_{j,k}$. Also, let $s_d^2=\sum_{k=1}^d\sigma_k^2$, and denote by $F_d$ the CDF of the normalized sum $\frac{1}{s_d}\sum_{k=1}^dA_k$. Then,
\begin{equation*}\label{eq11}
    F_d(x)=\Phi(x)-\phi(x)\Bigg[\frac{\lambda_3H_2(x)}{6s_d^3}
    +\frac{\lambda_4H_3(x)}{24s_d^4}+\frac{\lambda_3^2H_5(x)}{72s_d^6}+...\Bigg]
\end{equation*}
for all $x\in\mathbb{R}$, where $\Phi$ is the normal CDF, $\phi(x)=\frac{1}{\sqrt{2\pi}}e^{-\frac{1}{2}x^2}$ is the standard normal density, and $H_d$ is the $d$-th Hermite polynomial.
\end{theorem}

\section{Probabilistic Estimates}
Let us introduce Fristedt's probabilistic model for random partitions [6]. According to this model, the number of parts in a random partition equal to some integer $k\in\mathbb{N}$, denoted by $Z_k$, are independent random variables with geometric distribution. In particular, $Z_k\sim\mathrm{Geom}(1-q^k)$, for all $k\in\mathbb{N}$ and some $q\in(0,1)$. \\
In our work, we consider partitions $\lambda\in\mathcal{P}_\mu$, for some $\mu\in M$. Then, the condition $|\lambda|=n$ translates to
\begin{equation*}
    n=\sum_{k\in\mu(\mathbb{N})}k\mathbb{E}[Z_k]=\sum_{m=1}^{\infty}\mu(m)\mathbb{E}\big[Z_{\mu(m)}\big]=\sum_{m=1}^\infty\frac{\mu(m)q^{\mu(m)}}{1-q^{\mu(m)}}
\end{equation*}
where we have recalled that $\mathbb{E}(Z_k)=\frac{q^k}{1-q^k}$. It is convenient to set $q=\exp(-\alpha)$, then we obtain
\begin{equation}\label{eq1}
    n=\displaystyle\sum_{m=1}^\infty\frac{\mu(m)}{e^{\alpha\mu(m)}-1}.
\end{equation}
It is easy to see that equation~\eqref{eq1} has a unique solution $\alpha(n,\mu)>0$ for all $n\in\mathbb{N}, \mu\in M$. This follows from the Intermediate Value Theorem and the strict monotonicity of the r.h.s. on $\alpha$. Furthermore, we have the limit $\alpha\to0$ as $n\to\infty$, for all $\mu\in M$. \\
Additionally, notice that every $\mu\in M$ is an invertible function satisfying $\mu(x)=\Omega(x)$, this can be shown by the strict monotonicity and integer-mapping constraints in definition~\ref{def1}. Then, we have the following Lemma.
\begin{lemma}\label{lemma1}
    Let $\alpha$ be the solution of equation~\eqref{eq1} for some $n\in\mathbb{N},\mu\in M$, and let $y=y(n)$ satisfy $\alpha y\to\infty$ as $n\to\infty$. Then the largest part $Y_1$ of a random partition from $\mathcal{P}_\mu$ satisfies
    \begin{equation}\label{eq2}
        \mathbb{P}\big(\eta(Y_1)\le y\big)\xrightarrow{n\to\infty}e^{-e^{-y}}
    \end{equation}
    where $\eta(y)=\alpha y + \log\big(\alpha\mu'(\mu^{-1}(y))\big)$.
\end{lemma}
\begin{proof}
Let $\lambda\in\mathcal{P}_\mu$ and denote by $Y_1$ its largest part. Notice that $Y_1=\max\{k:Z_k>0\}$. According to Frisdedt's model, $Z_k\sim\mathrm{Geom}(1-q^k)$ are independent. Then $\mathbb{P}(Z_k>0)=q^k$ for all $k\in\mathbb{N}$. Thus
\begin{equation*}
    \mathbb{P}(Y_1\le y)=\mathbb{P}(Z_k=0, \ \forall k>y)=\prod_{\substack{m\in\mathbb{N} \\ \mu(m)>y}}\big(1-q^{\mu(m)}\big)
\end{equation*}
Setting $q=\exp(-\alpha)$, and recalling that $\log(1-x)=-x+O(x^2)$ for $x\to0$, we deduce
\begin{equation*}
\begin{split}
    \log\mathbb{P}(Y_1\le y)&=\sum_{\substack{m\in\mathbb{N} \\ \mu(m)>y}}\log\big(1-e^{-\alpha\mu(m)}\big)=-\sum_{\substack{m\in\mathbb{N} \\ \mu(m)>y}}e^{-\alpha\mu(m)} + \sum_{\substack{m\in\mathbb{N} \\ \mu(m)>y}}O\big(e^{-2\alpha\mu(m)}\big).
\end{split}
\end{equation*}
The error sum is of order $O(e^{-2\alpha y}/\alpha)$, and is therefore negligible compared to the main term. \\
Denote by $x_0=\mu^{-1}(y)$, where $\mu^{-1}$ is the inverse function of $\mu$. Then the main sum satisfies
\begin{equation*}
    \sum_{\substack{m\in\mathbb{N} \\ \mu(m)>y}}e^{-\alpha\mu(m)}=\int_{x_0}^\infty e^{-\alpha\mu(x)}dx+O(e^{-\alpha y})
\end{equation*}
A first-order Taylor expansion gives $\mu(x)=y+\mu'(x_0)(x-x_0)+O((x-x_0)^2)$. Hence
\begin{equation*}
    \int_{x_0}^\infty e^{-\alpha\mu(x)}dx=e^{-\alpha y}\int_{x_0}^\infty e^{-\alpha\mu'(x_0)(x-x_0)+O(\alpha(x-x_0)^2)}dx=\frac{e^{-\alpha y}}{\mu'(x_0)}(1+O(\alpha))
\end{equation*}
Combining the above estimates yield
\begin{equation*}
    \log\mathbb{P}(Y_1\le y)=-\frac{e^{-\alpha y}}{\mu'(\mu^{-1}(y))}(1+O(\alpha)).
\end{equation*}
Thus, by setting $\eta(y)=\alpha y + \log\big(\alpha\mu'(\mu^{-1}(y))\big)$, the rhs becomes $-e^{-\eta(y)}\big(1+O(\alpha)\big)$. Since $\eta$ is an invertible function for all $\mu\in M$ and $\alpha>0$, and since $\alpha\xrightarrow{n\to\infty}0^+$, we are done. 
\end{proof}

\begin{lemma}\label{lemma2}
    Let $\alpha$ be the solution of equation~\eqref{eq1} for some $n\in\mathbb{N},\mu\in M$, and let $x=x(n)$ satisfy $\alpha x\to\infty$ as $n\to\infty$. Then the number of parts $X_1$ of a random partition from $\mathcal{P}_\mu$ satisfies
    \begin{equation}\label{eq3}
        \mathbb{P}(\eta\circ\mu(X_1)\le x\big)\xrightarrow{n\to\infty}e^{-e^{-x}}
    \end{equation}
    where $\eta(x)=\alpha x+\log\big(\alpha\mu'(\mu^{-1}(x))\big)$.
\end{lemma}
\begin{proof}
    In order to prove the lemma, we shall denote by $P_\mu(n,x)$ the number of partitions of $n$ from $\mathcal{P}_\mu$ with at most $x$ parts, for all $n,x\in\mathbb{N}, \mu\in M$. In a similar manner as in the work of Szekeres [7], we obtain the generating function
    \begin{equation}\label{eq4}
        P_\mu(n,x)=\frac{1}{2\pi i}\oint_Cz^{-n-1}\prod_{m=x+1}^\infty\frac{1-z^{\mu(m)+\mu(x)}}{1-z^{\mu(m)}}dz
    \end{equation}
    where $C$ is a circular contour on the complex plane that centers at the origin and has a radius of $q=\exp(-\alpha)$. This is a trivial generalization of the work of Almkvist and Andrews [8], which considered the case $\mu(m)\equiv m$. Using the saddle-point method, used in [9], it is enough to evaluate the logarithm of the integrand and its derivatives in order to calculate the integral. After substituting $z=e^{-\alpha+i\theta}$, we may denote the logarithm of the integrand by $G(\theta,n,x)$ and obtain the following
\begin{equation*}
    G(\theta,n,x)=(\alpha-i\theta)n-\sum_{m=x+1}^\infty\log\Big(1-e^{(-\alpha+i\theta) \mu(m)}\Big)+\sum_{m=x+1}^\infty\log\Big(1-e^{(-\alpha+i\theta)(\mu(m)+\mu(x))}\Big).
\end{equation*}
Since we consider $x(n)$ for which $\alpha x\to\infty$ as $n\to\infty$, the second summation is negligible compared to the first summation, in the limit of large $n$. Furthermore, notice that for $\theta=0$, the first summation is of the same form as the sum we evaluated in the proof of Lemma~\ref{lemma1}, thus, we obtain in a similar manner
\begin{equation*}
    G(0,n,x)\xrightarrow{n\to\infty}\alpha n - \frac{e^{-\alpha\mu(x)}}{\alpha\mu(x)}(1+O(\alpha)).
\end{equation*}
Furthermore, it can be seen that the dependence on $x$ of the second derivative $\frac{d^2}{d\theta^2}G(0,n,x)$ is of order $O(e^{-\alpha x})$, and thus negligible. Richmond [10] has performed a complete derivation of this part for the case $\mu(x)\equiv x$. Since the derivation for general $\mu$ remains similar, we shall skip it. Then, according to the saddle-point theorem, we find with a good approximation the dependence of $P_\mu(n,x)$ on $x$, for large $n$:
\begin{equation*}
    \log P_\mu(n,x)\propto -\frac{e^{-\alpha\mu(x)}}{\alpha\mu(x)}(1+O(\alpha)).
\end{equation*}
We recall that $P_\mu(n,x)$ is the number of partitions of $n$ from $\mathcal{P}_\mu$ with at most $x$ parts, then the probability that a partition has at most $x$ parts is $\mathbb{P}(X_1\le x)\propto P_\mu(n,x)$. And thus we conclude that the variable $\eta\circ\mu(X_1)$ approaches a Gumbell distribution in the limit of large $n$, thus we are done.
\end{proof}
\begin{proposition}\label{prop1}
    Let $\mu\in M$ and let $n,k\in\mathbb{N}$. Indicate by $X_k$ the number of parts at least $k$ in a random partition of $n$ from the set $\mathcal{P}_\mu$, and denote by $Y_k$ the $k$-th largest part in the partition. Then the density distributions of $\mu(X_k)$ and $Y_k$ become equal as $n\to\infty$, and tend to
    \begin{equation*}\label{prop1}
        \frac{1}{(k-1)!}e^{-e^{-x}}e^{-kx}, \ \forall x\in\mathbb{R}, k\in\mathbb{N}.
    \end{equation*}
\end{proposition}
\begin{proof}
    From Lemma~\ref{lemma1}, the random variable $\eta(Y_1)$ approaches a Gumbell distribution in the limit of large $n$. Then from Gumbell order statistics, we deduce that the $k$-th largest part in a random partition from $\mathcal{P}_\mu$ satisfies
\begin{equation}\label{eq5}
    \mathbb{P}\big(\eta(Y_k)\le y\big)\xrightarrow{n\to\infty}\frac{1}{(k-1)!}\int_{-\infty}^ye^{-e^{-u}}e^{-ku}du
\end{equation}
for all $k\in\mathbb{N}$, and $y=y(n)$ satisfying $\alpha y\xrightarrow{n\to\infty}\infty$. \\
In addition, according to Lemma~\ref{lemma2}, the random variable $\eta\circ\mu(X_1)$ approaches a Gumbell distribution in the limit of large $n$. Then, in the same manner, we find
\begin{equation}\label{eq6}
    \mathbb{P}\big(\eta\circ\mu(X_k)\le x\big)\xrightarrow{n\to\infty}\frac{1}{(k-1)!}\int_{-\infty}^xe^{-e^{-u}}e^{-ku}du
\end{equation}
for all $k\in\mathbb{N}$, and $x=x(n)$ satisfying $\alpha x\xrightarrow{n\to\infty}\infty$. \\
Thus, the variables $\mu(X_k), Y_k$ have a similar asymptotic distribution, and their density distribution approaches
\begin{equation*}
    \frac{1}{(k-1)!}e^{-e^{-x}}e^{-kx}.
\end{equation*}
\end{proof}

\section{Restricted Graphical Partitions}
In the previous section, we studied the distributions of $X_k,Y_k$ of random partitions with restricted parts. In this section, we enumerate the graphical partitions of an integer $n$, under restrictions placed on their parts. To our knowledge, it has not been done before. Firstly, we shall introduce the following notation
\begin{notation}
    Let $n\in\mathbb{N},\mu\in M$, and let $\alpha$ be the solution of equation~\eqref{eq1}. Let $X$ be a random variable with density distribution function, for some $k\in\mathbb{N}$,
    \begin{equation*}
        \frac{1}{(k-1)!}e^{-e^{-x}}e^{-kx}, \ \forall x\in\mathbb{R}
    \end{equation*}
    Then for all $m\in\mathbb{N}$, denote
    \begin{equation}\label{eq7}
        \Delta_{k,m}=\mathbb{E}\Big[\big(\eta^{-1}(X)-\mathbb{E}[\eta^{-1}(X)]\big)^m\Big]-\mathbb{E}\Big[\big(\mu^{-1}\circ\eta^{-1}(X)-\mathbb{E}[\mu^{-1}\circ\eta^{-1}(X)]\big)^m\Big],
    \end{equation}
    where $\eta(x)=\alpha x+\log\big(\mu'(\mu^{-1}(x))\big)$.
\end{notation}
\begin{lemma}\label{lemma3}
    Let $\mu\in M$ be such that $\log\mu'(x)=o(x)$, and let $n\in\mathbb{N}$, denote by $\alpha$ the solution of equation~\eqref{eq1}. Then for $m=2,3,4$ and sufficiently large $k$
    \begin{equation}\label{eq8}
        \Delta_{k,m}\xrightarrow{n\to\infty}\frac{C_m}{\alpha^m k^{m-1}}
    \end{equation}
    where $C_m$ are constants independent of $k$.
\end{lemma}
\begin{proof}
    Since $\mu\in M$ we have $\mu^{-1}(x)=O(x)$, so under the assumption of $\log\mu'(x)=o(x)$ we obtain $\log\big(\mu'(\mu^{-1}(x))\big)=o(x)$. Thus, $\eta(x)\approx \alpha x$ for sufficiently large $x$, with error $o(x)$. Then, the inverse function satisfies $\eta^{-1}(x)\approx\alpha^{-1}x$. \\
    Denote by $X$ a random variable with density distribution function $f_k(x)$, for some $k\in\mathbb{N}$. Then we find for any $k\in\mathbb{N}$ and $m=1,2,3$, by equation~\eqref{eq7}
    \begin{equation*}
        \Delta_{k,m}\approx\frac{1}{\alpha^m}\mathbb{E}\Big[\big(X-\mathbb{E}[X]\big)^m\Big]-\mathbb{E}\Bigg[\bigg(\mu^{-1}\Big(\frac{1}{\alpha}X\Big)-\mathbb{E}\Big[\mu^{-1}\Big(\frac{1}{\alpha}X\Big)\Big]\bigg)^m\Bigg]
    \end{equation*}
    where the error decreases as $n\to\infty$. Notice that for large $n$ we have $\mu^{-1}(\alpha^{-1}X)=O(\alpha^{-1})$, then overall $\Delta_{k,m}=O(\alpha^{-m})$ as $n\to\infty$. Furthermore, if $\mu^{-1}(x)=o(x)$ then the second expected value is negligible compared to the first, for sufficiently large $n$, and if $\mu^{-1}(x)=\Theta(x)$ then it approaches a constant multiple of the first expected value. Then in general, the order of growth of $\Delta_{k,m}$ for large $n$ is determined only by the first expected value. Notice that
    \begin{equation*}
        \mathbb{E}[X]=\frac{1}{(k-1)!}\int_{-\infty}^\infty xe^{-e^{-x}}e^{-kx}dx=-\psi(k)\approx-\log k
    \end{equation*}
    where $\psi$ is the digamma function and the last step is an approximation for large $k$. Thus, one can check that
    \begin{equation*}
        \mathbb{E}\Big[\big(X-\mathbb{E}[X]\big)^2\Big]=\frac{1}{(k-1)!}\int_{-\infty}^\infty\big(x+\psi(k)\big)^2e^{-e^{-x}}e^{-kx}dx=\psi'(k)\approx\frac{1}{k}.
    \end{equation*}
    and in the same manner, $\mathbb{E}\Big[\big(X-\mathbb{E}[X]\big)^m\Big]\approx\frac{1}{k^{m-1}}$ also for $m=3,4$. This concludes the proof.
\end{proof}
\begin{lemma}\label{lemma4}
    Let $k,n\in\mathbb{N}, \mu\in M$. Also, let $R_k$ be the $k$-th successive rank of a random partition of $n$, from $\mathcal{P}_\mu$. Then $R_k$ is a random variable that satisfies for all $1\le m\le4$
    \begin{equation}\label{eq9}
        \mathbb{E}\Big[\big(R_k-\mathbb{E}[R_k]\big)^m\Big]=\Delta_{k,m}.
    \end{equation}
\end{lemma}
\begin{proof}
    For $m=1$ the proof is trivial, since equation~\eqref{eq7} nullifies in that case. \\
    As mentioned above, the $k$-th rank of a partition satisfies $R_k=Y_k-X_k$. For a random partition from $\mathcal{P}_\mu$, we denote by $\mathcal{R}_k(r)$ the density distribution of $R_k$, for all $k\in\mathbb{N}$. Then, from Theorem~\ref{theo1} we deduce
    \begin{equation*}
    \mathcal{R}_k(r)=\int_{-\infty}^\infty f_k\big(\eta(r+x)\big)f_k\big(\eta\circ\mu(x)\big)J(r,x)dx
    \end{equation*}
    where the density distributions are
    \begin{equation*}
        f_k(x)\equiv\frac{1}{(k-1)!}e^{-e^{-x}}e^{-kx}, \ \forall k\in\mathbb{N}
    \end{equation*}
    and the Jacobian is $J(r,x)=\frac{d}{dx}\eta(r+x)\frac{d}{dx}\eta\big(\mu(x)\big)$. Then, substituting $\eta(r+x)\mapsto u$, we obtain the following
    \begin{equation*}
        \mathbb{E}[R_k]=\int_{-\infty}^\infty r\mathcal{R}_k(r)dr = \int_{-\infty}^\infty dx\Bigg[f_k\big(\eta\circ\mu(x)\big)\frac{d}{dx}\eta\circ\mu(x)\int_{-\infty}^\infty du(\eta^{-1}(u)-x)f_k(u)\Bigg].
    \end{equation*}
    Setting $\eta\circ\mu(x)\mapsto t$ and simplifying, we find
    \begin{equation*}
        \mathbb{E}[R_k]=\int_{-\infty}^\infty\eta^{-1}(u)f_k(u)du-\int_{-\infty}^\infty\mu^{-1}\circ\eta^{-1}(t)f_k(t)dt.
    \end{equation*}
    Then, by conducting a similar calculation for $m=2,3,4$, one can check that
    \begin{equation}\label{eq10}
        \mathbb{E}\big[R_k^m\big]=\int_{-\infty}^\infty\Big(\eta^{-1}(u)-\mu^{-1}\circ\eta^{-1}(u)\Big)^mf_k(u)du , \ \ 1\le m\le4.
    \end{equation}
    From this, the rest of the proof is trivial from the definition of $\Delta_{k,m}$.
\end{proof}
Now, we shall prove Theorem~\ref{theo1}.
\begin{theorem}\label{theo1}
    Let $\mu\in M$ and let $n\in\mathbb{N}$ be an even integer. Consider a partition of $n$ with Durfee square side length $D$, chosen uniformly at random from $\mathcal{P}_\mu$. Then, the probability that this partition is graphical is bounded above by
    \begin{equation}
        \frac{C}{\log^{3/2}D}
    \end{equation}
    for sufficiently large $D$, with $C=\sqrt{\frac{32}{\pi}}e^{-5/4}\approx 0.9144$.
\end{theorem}
\begin{proof}
    Denote $A_k=R_k-\mathbb{E}[R_k]$ for any integer $k\in\mathbb{N}$, then according to Lemma~\ref{lemma4}
    \begin{equation*}
        \mathbb{E}[A_k]=0, \ \ \mathbb{E}\big[A_k^2\big]=\Delta_{k,2}, \ \ \mathbb{E}\big[A_k^3\big]=\Delta_{k,3}.
    \end{equation*}
    Let $\alpha=\alpha(d,\mu)$ be the solution of equation~\eqref{eq1} for the integer $d$ and function $\mu\in M$. Then from Lemma~\ref{lemma3}, using the same notation used in Theorem~\ref{Edge}, we evaluate
    \begin{equation*}
    \begin{split}
        &s_d^2\equiv\sum_{k=1}^d\Delta_{k,2}\xrightarrow{d\to\infty}\sum_{k=1}^d\frac{C}{\alpha^2k}=\Theta\big(\alpha^{-2}\log d\big), \\ &\lambda_m\equiv\sum_{k=1}^d\Delta_{k,m}\xrightarrow{d\to\infty}\sum_{k=1}^d\frac{C_m}{\alpha^mk^{m-1}}=\Theta(\alpha^{-m}),
    \end{split}
    \end{equation*}
    for $m=3,4$, where $C,C_m$ are constants. \\
    Let $F_d$ be the CDF of the normalized sum $\frac{1}{s_d}\sum_{k=1}^dA_k$, for some $d\in\mathbb{N}$, then from Proposition~\ref{prop1} we find for all $x\in\mathbb{R}$
    \begin{equation*}
        F_d(x)=\Phi(x)-\frac{\phi(x)H_2(x)}{(\log d)^{3/2}}-\phi(x)p(x)\cdot O\bigg(\frac{1}{\log^2d}\bigg)
    \end{equation*}
    where $p(x)$ is a polynomial. Since $\phi(x)$ is a gaussian, the rhs is bounded from above, thus, for sufficiently large $d$
    \begin{equation*}
        \sup_x\big|F_d(x)-\Phi(x)\big|=\Theta\Big((\log d)^{-3/2}\Big).
    \end{equation*}
    Where the underlying constant $C$ is:
    \begin{equation}
        C=\sup_x\phi(x)H_2(x)=\sup_x\frac{1}{\sqrt{2\pi}}e^{-x^2/2}(4x^2-2)=\sqrt{\frac{32}{\pi}}e^{-5/4}\approx0.9144.
    \end{equation}
    Then, we obtain
    \begin{equation*}
    \begin{split}
        &\mathbb{P}\bigg[\sum_{k=1}^dR_k\le-d\bigg]=\mathbb{P}\bigg[\frac{1}{s_d}\sum_{k=1}^dA_k\le-\frac{d}{s_d}-\frac{1}{s_d}\sum_{k=1}^d\mathbb{E}[R_k]\bigg] \\
        &=F_d\bigg(-\frac{d}{s_d}-\frac{1}{s_d}\sum_{k=1}^d\mathbb{E}[R_k]\bigg)\le \Phi\bigg(-\frac{d}{s_d}-\frac{1}{s_d}\sum_{k=1}^d\mathbb{E}[R_k]\bigg)+C(\log d)^{-3/2}.
    \end{split}
    \end{equation*}
    Considering equation~\eqref{eq10} with $m=1$, we see that $\mathbb{E}[R_k]\ge0, \forall k\in\mathbb{N},\mu\in M$, then the $\Phi$ term becomes negligible for large $d$, and thus
    \begin{equation*}
        \mathbb{P}\bigg[\sum_{k=1}^dR_k\le-d\bigg]\le C(\log d)^{-3/2}.
    \end{equation*}
    Hence, from Theorem~\ref{Nash}, the probability that a partition from $\mathcal{P}_\mu$ with Durfee square $D$ is graphical satisfies
    \begin{equation*}
        \mathbb{P}\bigg[\sum_{k=1}^dR_k\le-d, \ \forall d\le D\bigg]\le C\big(\log D\big)^{-3/2},
    \end{equation*}
    for sufficiently large $D$. \\
   Thus we are done.
\end{proof}
We conclude that the probability of a random integer partition being graphical is asymptotically bounded by its Durfee square side length. For unrestricted partitions, the typical Durfee square side length is $\sqrt{6n}\log2/\pi$ [11], however, it remains unknown what is the typical value for general $\mu\in M$. In the remainder of this article we derive the typical side length of the Durfee square when general restrictions are placed on the parts. Specifically, we prove the following theorem.
\begin{theorem}\label{theo2}
    Let $\mu\in M$ and let $n$ be an even integer. Then the most-probable Durfee square side length $D$ of a random partition of $n$ from $\mathcal{P}_\mu$ is the solution of the trancedental equation
    \begin{equation}\label{D-eq}
        \beta D+\frac{\log(1-e^{-\beta D})}{\mu'(\mu^{-1}(D))}=0.
    \end{equation} \\
    And the variance of $D$ is $O(\beta^3)$, where $\beta>0$ is defined by
    \begin{equation}\label{beta}
        \frac{n-D^2}{2}=\sum_{\substack{m\in\mathbb{N} \\ \mu(m)\le D}}\frac{\mu(m)}{e^{\beta\mu(m)}-1}.
    \end{equation}
\end{theorem}

\vspace{0.5cm}

In order to prove Theorem~\ref{theo2}, we shall show the following lemma first.
\begin{lemma}\label{D-lemma}
    Denote by $N(n,D,\mu)$ the number of partitions of $n$ from $\mathcal{P_{\mu}}$ with Durfee square side length $D$. Then
    \begin{equation}
        N(n,D,\mu)=e^{\beta(n-D^2)}\prod_{\substack{m\in\mathbb{N} \\ \mu(m)\le D}}\big(1-e^{-\beta\mu(m)}\big)^{-2} \cdot O(\beta^{3/2}),
    \end{equation}
    where $\beta$ is defined by~\eqref{beta}.
\end{lemma}

\begin{proof}
    Firstly, we denote by $N(n,D,\mu)$ the number of partitions of $n$ from $\mathcal{P_{\mu}}$ with a Durfee square side length of $D$. Andrews [12] showed that for unrestricted partitions, with $\mu(x)\equiv x$, it holds that
    \begin{equation*}
        N(n,D)=\big[q^{n-D^2}\big]\prod_{\substack{m=1}}^\infty(1-q^{m})^{-2},
    \end{equation*}
    where $q\in(0,1)$. Then by a natural generalization, for any $\mu\in M$ the following holds
    \begin{equation}
        N(n,D,\mu)=\big[q^{n-D^2}\big]\prod_{\substack{m\in\mathbb{N} \\ \mu(m)\le D}}\big(1-q^{\mu(m)}\big)^{-2}.
    \end{equation}
    Thus, we obtain
    \begin{equation*}
        N(n,D,\mu)=\frac{1}{2\pi i}\oint q^{-n+D^2-1}\prod_{\substack{m\in\mathbb{N} \\ \mu(m)\le D}}\big(1-q^{\mu(m)}\big)^{-2}dq
    \end{equation*}
    where the integration is taken along an infinitesimally small circular contour centered at the origin. Setting $q=e^{-\beta}$ we have
    \begin{equation*}
        N(n,D,\mu)=\frac{1}{2\pi i}\oint \exp\bigg[{-\beta(n-D^2)-\sum_{\substack{m\in\mathbb{N} \\ \mu(m)\le D}} 2\log\big(1-e^{-\beta \mu(m)}\big)}\bigg]d\beta.
    \end{equation*}
    This intgeral can be approximated using the saddle point method, in a similar manner as conducted above. We shall denote the logarithm of the integrand by $\Phi(\beta)$, and find its maximum point by differentiating:
    \begin{equation*}
        \Phi'(\beta)=n-D^2-\sum_{\substack{m\in\mathbb{N} \\ \mu(m)\le D}}\frac{2\mu(m)}{e^{\beta\mu(m)}-1}=0.
    \end{equation*}
    Hence, we obtained an equation for the parameter $\beta$, in terms of $n, D$ and $\mu$. In a similar behavior as for the parameter $\alpha$, we have $\beta\to0$ as $(n-D^2)\to\infty$, for any $\mu\in M$. \\
    \\
    Furthermore, the second derivative satisfies
    \begin{equation}\label{second-der}
        \Phi''(\beta)=\sum_{\substack{m\in\mathbb{N} \\ \mu(m)\le D}}\frac{2\mu^2(m)e^{\beta\mu(m)}}{(e^{\beta\mu(m)}-1)^2}=\int_0^D\frac{2x^2e^{\beta x}}{(e^{\beta x} -1)^2}\frac{dx}{\mu'(\mu^{-1}(x))} + O(\beta)=O\Big(\frac{1}{\beta^3}\Big).
    \end{equation}
    Thus, according to the saddle point method,
\begin{equation*}
    N(n,D,\mu)=\frac{e^{\Phi(\beta)}}{\sqrt{\Phi''(\beta)}}\cdot O(1)=e^{\beta(n-D^2)}\prod_{\substack{m\in\mathbb{N} \\ \mu(m)\le D}}\big(1-e^{-\beta\mu(m)}\big)^{-2} \cdot O(\beta^{3/2}),
\end{equation*}
and the proof is complete.
\end{proof}
Now we shall prove Theorem~\ref{theo2}:
\begin{proof}
    The probability that a random partition of $n$ from $\mathcal{P}_\mu$ has a Durfee square side length $D$ is
\begin{equation*}
    \mathbb{P}[D|n,\mu]=\frac{N(n,D,\mu)}{|\mathcal{P}_\mu(n)|}.
\end{equation*}
Notice that the denominator does not depend on $D$, then, in order to find the value of $D$ with maximal probability, we shall consider only the numerator. From lemma~\ref{D-lemma} we have, by approximating the summation by an integral
\begin{equation}
    \log(N(n,D,\mu))=\beta(n-D^2)-2\int_0^{\mu^{-1}(D)}\log(1-e^{-\beta\mu(x)})dx+O(\log\beta).
\end{equation}
Although $D$ is a discrete integer, we shall treat it as a continuous parameter and differentiate by it to find the maximum, that is, we demand:
\begin{equation*}
    \beta'(n-D^2)-2\beta D-\frac{2\log\big(1-e^{-\beta D}\big)}{\mu'(\mu^{-1}(D))}-2\beta'\int_0^{\mu^{-1}(D)}\frac{\mu(x)}{e^{\beta\mu(x)}-1}dx+O\Big(\frac{\beta'}{\beta}\Big)=0.
\end{equation*}
Where $\beta'\equiv\partial_D\beta$. \\
Notice that the first and last terms cancel out, up to a negligible error of $O(\beta')$. Furthermore, the error $O(\beta'/\beta)$ is negligible compared to the second and third terms, this can be seen from~\eqref{beta}. Thus, we conclude that the side length $D$ with maximum probability satisfies:
\begin{equation*}
    \beta D+\frac{\log\big(1-e^{-\beta D}\big)}{\mu'(\mu^{-1}(D))}=0.
\end{equation*}
Furthermore, the variance of $D$ is simply $1/\Phi''(\beta)$, which we concluded in~\eqref{second-der} is $O(\beta^3)$. \\
Thus, we are done.
\end{proof}

\vspace{0.75cm}

To emphasize the consequences of the results shown in this article, we shall consider some specific cases. Notice that for $\mu(x)\equiv ax$, with $a\ge1$, equation~\eqref{D-eq} becomes
\begin{equation*}
    \beta D + \frac{1}{a}\log(1-e^{-\beta D})=0,
\end{equation*}
where $\beta=\frac{\pi}{a\sqrt{6(n-D^2)}}\approx\frac{\pi}{a\sqrt{6n}}$ with negligible error. Setting $y=e^{-\beta D}$ yields
\begin{equation}
    y^a+y-1=0.
\end{equation}
In particular, for $a=1$ we have $y=\frac{1}{2}$, and we obtain the well-known typical value $D=\sqrt{6n}\log2/\pi$. \\
\\
Another nice example is for $a=2$, that is, when the parts are restricted to even integers, in this case we obtain $D=\sqrt{24n}\log(\varphi)/\pi$, where $\varphi$ is the golden ratio. \\
\\
Then, using Theorem~\ref{theo1} and the results above, we can generally conclude the following theorem~\ref{theo3}, regrading restrictions of the form $\mu(x)\equiv ax$, for any $a>1$.

\newpage

\begin{theorem}\label{theo3}
    The probability that a random partition of $n$ is graphical, given that its parts are taken from a set whose elements grow asymptotically linearly, is bounded from above by
\begin{equation}
    \frac{C_\mathrm{lin}}{\log^{3/2}n}\Big(1+O\Big(\frac{1}{\log n}\Big)\Big)
\end{equation}
where $C_\mathrm{lin}=\sqrt{\frac{8}{\pi}}e^{-5/4}\approx0.4572$.
\end{theorem}

Now consider partitions whose parts are restricted to perfect squares, that is, $\mu(x)\equiv x^2$. In this case, we have $\beta=(4n/\sqrt{\pi}\zeta(3/2))^{-2/3}$ for sufficiently large $n$, where $\zeta$ is the Riemann zeta function. Furthermore, equation~\eqref{D-eq} becomes
\begin{equation*}
    2\beta D^{3/2} + \log(1-e^{-\beta D})=0.
\end{equation*}
An asymptotic solution yields, with a negligible relative error
\begin{equation*}
    D=\Big(\frac{4n}{\sqrt{3\pi}\zeta(3/2)}\Big)^{4/9}\log^{2/3}n.
\end{equation*}
Thus, we conclude
\begin{theorem}\label{theo4}
    The probability of a random partition of $n$ with parts restricted to perfect squares to be graphical is bounded from above by
    \begin{equation}
    \frac{C_\mathrm{sq}}{\log^{3/2}n}\Big(1+O\Big(\frac{\log\log n}{\log n}\Big)\Big).
\end{equation}
where $C_\mathrm{sq}=\sqrt{\frac{729}{2\pi}}e^{-5/4}\approx3.086$.
\end{theorem}

\vspace{0.5cm}

\section*{Conclusions}
In this work, we established a universal upper bound for the probability that a random integer partition is graphical. This bound depends solely on the partition's Durfee square side length and is invariant with respect to $n$ and any restrictions on its parts. Because the side length of the Durfee square itself depends on these criteria, we derived its typical size under general restrictions on the parts. These results enabled us to formulate bounds on the probability of graphicality that depend directly on $n$ and the given restrictions. \\
Looking ahead, it would be worthwhile to establish tighter bounds on this probability and to extend our framework to directed graphs and multigraphs.

\end{document}